\documentclass[12pt]{article}

\newtheorem{proposition}{Proposition}
\newtheorem{theorem}{Theorem}
\newtheorem{corollary}{Corollary}
\newtheorem{lemma}{Lemma}
\newtheorem{remark}{Remark}

\newtheorem{example}{Example}
\usepackage{graphics}

\usepackage{color}
\usepackage[%
breaklinks%
,colorlinks%
,linkcolor=red
,anchorcolor=blue%
,pagecolor=blue%
,citecolor=blue%
,bookmarks=ture%
]{hyperref}

\def\qed{\hfill  \framebox(5,5){}}
\def\deg{{\rm deg}}
\def\cP{{\mathcal P}}
\def\cQ{{\mathcal Q}}

\expandafter\chardef\csname pre amssym.def
at\endcsname=\the\catcode`\@ \catcode`\@=11

\def\undefine#1{\let#1\undefined}
\def\newsymbol#1#2#3#4#5{\let\next@\relax
 \ifnum#2=\@ne\let\next@\msafam@\else
 \ifnum#2=\tw@\let\next@\msbfam@\fi\fi
 \mathchardef#1="#3\next@#4#5}
\def\mathhexbox@#1#2#3{\relax
 \ifmmode\mathpalette{}{\m@th\mathchar"#1#2#3}%
 \else\leavevmode\hbox{$\m@th\mathchar"#1#2#3$}\fi}
\def\hexnumber@#1{\ifcase#1 0\or 1\or 2\or 3\or 4\or 5\or 6\or 7\or 8\or
 9\or A\or B\or C\or D\or E\or F\fi}

\font\tenmsa=msam10 \font\sevenmsa=msam7 \font\fivemsa=msam5
\newfam\msafam
\textfont\msafam=\tenmsa \scriptfont\msafam=\sevenmsa
\scriptscriptfont\msafam=\fivemsa \edef\msafam@{\hexnumber@\msafam}
\mathchardef\dabar@"0\msafam@39
\def\dashrightarrow{\mathrel{\dabar@\dabar@\mathchar"0\msafam@4B}}
\def\dashleftarrow{\mathrel{\mathchar"0\msafam@4C\dabar@\dabar@}}

\font\tenmsb=msbm10 \font\sevenmsb=msbm7 \font\fivemsb=msbm5
\newfam\msbfam
\textfont\msbfam=\tenmsb \scriptfont\msbfam=\sevenmsb
\scriptscriptfont\msbfam=\fivemsb \edef\msbfam@{\hexnumber@\msbfam}
\def\Bbb#1{\fam\msbfam\relax#1}

\def\deg{{\rm deg}}

\def\lcm{{\rm lcm}}

\def\cP{{\mathcal P}}

\def\ox{\,{\overline x}\,}

\def\para{\vspace{2.5 mm}}

\def\lcm{{\rm lcm}}

 \def\res{{\rm Res}}
\def\cP{{\mathcal P}}

\def\ox{\,{\overline x}\,}

\begin{document}

\title{Parametrization of  translational surfaces\thanks{This project is supported by Beijing Nova Program under Grant Z121104002512065.  The author S.P\'erez-D\'{\i}az is member of the Research Group ASYNACS (Ref. CCEE2011/R34)}}
\author{Sonia P\'erez-D\'{\i}az  \\
Dpto de F\'{\i}sica y Matem\'aticas \\
        Universidad de Alcal\'a \\
      E-28871 Madrid, Spain  \\
sonia.perez@uah.es
\and Li-Yong Shen\\
School of Mathematical Sciences  \\
      University of CAS \\
        Beijing, China \\lyshen@ucas.ac.cn
      }

\date{}          
\maketitle

\begin{abstract}
The algebraic translational surface is a typical modeling surface in
computer aided design and architecture industry. In this paper, we give a necessary and sufficient condition for that algebraic surface having a standard parametric representation and our proof is constructive. If the given algebraic surface is translational,  then we can  compute a standard parametric representation for the surface.
\end{abstract}



\section{Introduction}
In computer aided geometric design and geometric modeling, we often use some common surface to construct a complex model. These common  surfaces are generally called basic modeling surfaces, and they should have some advantages such as  simple and widely used.  The simple surfaces  refer the ones with low degree, for instance, quadratic surfaces \cite{ivan95,lazard02, wang}  and cubic surfaces \cite{bajaj98, berry01}. The widely used surface refer the ones common in industrial design, for instance,  ruled surfaces \cite{chen03,chen11,shen14}, swept surfaces \cite{schroeder94, weinert04} and translational surfaces \cite{glymph04, liu06}.  Then a primary task is to have a well study for these basic modeling surfaces.
 Certainly, how to represent these surface is the first problem in studying surface.

 \para

 As we know, two representation forms are usually
used as parametric form and   implicit form. For the parametric representation, it is a main popular
geometric representation in CG and CAD~\cite{handbook}. It is easy to render and is helpful for some geometric operations, such as the computation of curvature or bounds and the control of position or tangency. However, it is difficult in positional relationship determination and collision detection.  Another limitation is that the parametric surfaces have lower ability in topology comparing with the implicit surfaces. This is also the reason that modeling technology of implicit surfaces is increasing in more applications. However, it is difficult to classify and  control  the implicit surfaces. Hence, it is a problem to determining some basic modeling surfaces in the implicit form, and furthermore, find a parametric representation if it has.

 \para

In this paper, we prefer to study the translational surface which is commonly used in industrial design. The translational surface is a simple solution to the problem of interpolating a surface passing through two curves.   Hence, people generally give two space curves ${\cal C}_1$ and ${\cal C}_2$ with parametrizations ${\cal P}_1(t_1)$ and ${\cal P}_2(t_2)$, respectively, and the two curves define  a translational surface ${\cal P}(t_1,t_2)={\cal P}_1(t_1)+{\cal P}_2(t_2)$.
However, as we mentioned above, newly geometric modeling often begin with implicit equations \cite{turk02,yang05}. In this situation, for better  control and design,  it is necessary to find the basic modeling surface and compute their parametric representations if exist. The main task of this paper is to parametrize a given algebraic translational surface. Considering the parametric form of translational surface, the two variables are separated. Based on this fact we firstly find a space curve which can play the role ${\cal C}_1$.  Then we compute a parametrization ${\cal P}_1(t_1)$ for it. Successively, we find ${\cal C}_2$ and  compute a parametrization ${\cal P}_2(t_2)$ for it.  Finally, we get a parametrization  ${\cal P}(t_1,t_2)={\cal P}_1(t_1)+{\cal P}_2(t_2)$ for the given  algebraic translational surface.

 \para

The paper is organized as follows. First, we introduce some preliminaries and previous results  (see Section 2). In Section 3, we show a necessary and sufficient condition for that an algebraic surface has a parametric representation of the form ${\cal P}(t_1,t_2)={\cal P}_1(t_1)+{\cal P}_2(t_2)$. The proof is constructive, and then a method for computing $\cP$ is developed. Finally, we show how the computation of $\cP$ can be improved, computationally speaking. More precisely, the final problem consist in deciding the existence of two space curves and to compute a rational parametrization of them (see Section 4).


\section{Preliminaries}

Let $\cal V$ be a  surface over an algebraically closed field of characteristic zero $\Bbb K$, and let $f(\ox)\in {\Bbb K}[\ox],\,\,\ox:=(x_1, x_2, x_3)$, be the irreducible polynomial defining implicitly  $\cal V$.

\para

In the following, we analyze if $\cal V$ is a {\it translational surface}; i.e.   whether $\cal V$  admits a parametrization of the standard form
\begin{equation}\label{eq-form}
{\cal P}(t_1, t_2)={\cal P}_1(t_1)+{\cal P}_2(t_2),\end{equation} where
\[ {\cal P}_1(t_1)=(p_{11}(t_1), p_{12}(t_1), p_{13}(t_1))\in {\Bbb K}(t_1)^3\setminus{\Bbb K}^3,\]\[{\cal P}_2(t_2)=(p_{21}(t_2), p_{22}(t_2), p_{23}(t_2))\in {\Bbb K}(t_2)^3\setminus{\Bbb K}^3\]
and in the affirmative case, we compute it. We denote by ${\cal C}_i$ the space curve over $\Bbb K$ defined by the rational parametrization $\cP_i$, for $i=1,2$.

\para

Throughout this paper,  $\cal V$ is assumed  not to be a cylindrical surface. One can check that  $f(\ox)=0$ defines  a cylindrical surface if and only if there exists a constant vector $(a_1,a_2,a_3)\in {\Bbb K}^3\setminus\{(0,0,0)\}$ such that $\nabla f(\ox)\cdot (a_1, a_2, a_3)=  a_1f_{x_1}+a_2f_{x_2}+a_3f_{x_3}=0$, where  $f_{x_i}$ denotes the partial derivative of the polynomial $f$ w.r.t. the variable $x_i$. For this case, it is not difficult to get  a parametrization for $\cal V$ (see \cite{shen14b}). As a special case, we get that $\cal V$ is not a plane. The plane case  has a trivial solution since if $\cal V$ is defined by the polynomial $f(\ox)=m_1x_1+m_2x_2+m_3x_3+m_4,\,m_i\in {\Bbb K}$, a parametrization of the form given in Eq. (\ref{eq-form}) is  ${\cal P}(t_1, t_2)={\cal P}_1(t_1)+{\cal P}_2(t_2)=(t_1, t_2, -m_1t_1-m_2t_2-m_4)$, where ${\cal P}_1(t_1)=(t_1, 0, -m_1t_1-m_4),\,\,{\cal P}_2(t_2)=(0, t_2, -m_2t_2)$.

\para

In the following, we present some properties concerning the parametrizations $\cP_i(t_i)\in {\Bbb K}(t_i)^3\setminus{\Bbb K}^3$, and the space curves ${\cal C}_i,\,i=1,2$. These results,  will play an important role in  Sections 3 and 4.

\para

\begin{proposition}\label{P-properties}  If $\cal V$ is a translational surface, then the following properties hold:
\begin{itemize}
\item[1.] If $\cP$ is proper, then $\cP_1$ and $\cP_2$ are both proper parametrizations. 
\item[2.] There exist  $\cP_1$ and $\cP_2$  such that they are both proper parametrizations.
\item[3.] It holds that ${\cal C}_1\not={\cal C}_2$.
\item[4.] It holds that ${\cal C}_i$ is not a line, for $i=1,2$.
\end{itemize}
\end{proposition}

\vspace*{2mm}

\noindent {\bf Proof.}
 \begin{itemize}
\item[1.] Let us prove that $\cP_1$ is   a  proper parametrization. Similarly, one shows that $\cP_2$ is proper. Let us assume that $\cP_1$ is not proper.  Then, there exists $\phi_i(s_1)\in \overline{{\Bbb K}(s_1)},\,i=1,2,\,\phi_1\not=\phi_2$ such that $\cP_1(\phi_1(s_1))=\cP_1(\phi_2(s_1))=\cP_1(s_1)$ ($\overline{{\Bbb K}(s_1)}$ is the algebraic closure of ${\Bbb K}(s_1)$, and $s_1$ is a new variable). Thus,   $\cP(\phi_1(s_1),s_2)=\cP(\phi_2(s_1),s_2)=\cP(s_1,s_2)$ ($s_2$ is a new variable). This implies that $\cP$ is not proper. Therefore, we conclude that $\cP_1$ is proper which is a contradiction.
\item[2.] Let us assume that ${\cal P}_1(t_1)$ or ${\cal P}_2(t_2)$ is not proper. Then, we may reparametrize each of them using for instance the results in \cite{sonia06}. That is, there exists  ${\cal P}^*_i$  proper, and $R_i(t_i)\in {\Bbb K}(t_i)\setminus{\Bbb K}$ such that ${\cal P}^*_i(R(t_i))={\cal P}_i(t_i),\,i=1,2$. Under these conditions, we have that ${\cal P}^*(t_1, t_2)={\cal P}^*_1(t_1)+{\cal P}^*_2(t_2)$ satisfies that ${\cal P}^*(R_1(t_1), R_2(t_2))={\cal P}(t_1, t_2)$, and thus  ${\cal P}^*$ is a parametrization of $\cal V$ of the form given in Eq. (\ref{eq-form}).
\item[3.] Clearly ${\cal C}_1\not={\cal C}_2$ since $\cP$ parametrizes a surface $\cal V$.
\item[4.] Let us assume that ${\cal C}_2$ is a line (we reason similarly, if ${\cal C}_1$ is a line). Then, a proper parametrization of ${\cal C}_2$ is given by $\cP_2(t_2)=(a_1t_2+b_1,a_2t_2+b_2,a_3t_2+b_3)\in {\Bbb K}(t_2)^3\setminus{\Bbb K}^3$. Then, since $f(\cP(t_1, t_2))=0$, we get that $\nabla f(\cP(t_1, t_2))\,\cdot \,{\cal P}_2'(t_2)=\nabla f(\cP(t_1, t_2))\,\cdot \,(a_1, a_2, a_3)=0$. Thus, $\cP$ parametrizes the surface $\cal V$ defined by $f(\ox)$, and the surface defined by  $\nabla f(\ox)\,\cdot \,(a_1, a_2, a_3)$. Since $f$ is irreducible, we get that $\nabla f(\ox)\,\cdot \,(a_1, a_2, a_3)=0$ which is impossible since $\cal V$ is not  a cylindrical surface.
 \qed
\end{itemize}

\begin{remark}\label{R-properties} In the following, we compute $\cP_i,\,i=1,2$ being proper (see statement $2$ in Proposition \ref{P-properties}). However, we note that the properness of $\cP_i,\,i=1,2$ does not imply that $\cP$ is proper.
\end{remark}

\para

If  $\cal V$ is a translational surface, then $\cal V$    admits a parametrization of the form given in Eq. (\ref{eq-form}), and  ${\cal P}_2(t_2)\in {\Bbb K}(t_2)^3$ can be assumed to satisfied  some  additional properties.
In particular, we claim that there exists a $t_2^0\in \Bbb{K}$ such that  ${\cal P}_2(t_2^0)=(0,0,0)$, and  ${\cal P}'_2(t_2^0)\not=(0,0,0)$.
Without loss of generality, we prove the following lemma by setting $t_2^0=0$.

\para

\begin{lemma}\label{L-prop1}  Let  $\cal V$ be a translational surface. Then, there exists a proper parametrization   ${\cal P}_2(t_2)\in {\Bbb K}(t_2)^3$ such that ${\cal P}_2(0)=(0,0,0)$, and  ${\cal P}'_2(0)\not=(0,0,0)$.
\end{lemma}

\vspace*{2mm}

\noindent {\bf Proof.}  Let $\overline{{\cal P}}(t_1, t_2)=\widehat{\cal P}_1(t_1)+\overline{{\cal P}}_2(t_2)$ be a parametrization of $\cal V$ such that $\widehat{\cal P}_1(t_1), \overline{{\cal P}}_2(t_2)$ are proper (see statement $2$ in Proposition \ref{P-properties}). We consider a change of variable on the variable $t_2$ of the form $$\phi(t_2)=\frac{at_2+b}{ct_2+1}\in {\Bbb K}(t_2)\setminus{\Bbb K},\,\,\mbox{where}\,\,a-bc\not=0\quad \qquad \,\,\mbox{(I)}$$ and such that  $\overline{{\cal P}}_2$ is defined at $t_2=\phi(0)=b$. That is, $$b\not=r_i,\,\,i=1,\ldots,m,\,\, \mbox{and}\,\, p_2(r_i)=0,\,\,\mbox{where}\,\,p_2:=\lcm(p_{212},p_{222},p_{232}),\,p_{2j}=\frac{p_{2j1}}{p_{2j2}}\,\,\,\,\,  \mbox{(II)}.$$

Let  ${\cal P}(t_1,t_2)=\overline{{\cal P}}(t_1, \phi(t_2)) =\widehat{\cal P}_1(t_1)+ \widehat{\cal P}_2(t_2)$,   where $\widehat{\cal P}_2(t_2)=\overline{{\cal P}}_2(\phi(t_2))$. Note that  $\widehat{\cal P}_2(0)=\overline{{\cal P}}_2(b)$ ($\overline{{\cal P}}_2$ is defined at $b$),  and $\widehat{\cal P}_2$ is proper since  $\widehat{\cal P}_2(t_2)=\overline{{\cal P}}_2(\phi(t_2))$ and $\phi(t_2)$ is invertible and $\overline{{\cal P}}_2(t_2)$ is proper. Furthermore, ${\cal P}$ is a parametrization of $\cal V$   ($\overline{{\cal P}}$ is  a parametrization of $\cal V$, and ${\cal P}(t_1,t_2)=\overline{{\cal P}}(t_1, \phi(t_2))$). \para

\noindent
In addition,
\[\widehat{\cal P}'_2(0)=\overline{{\cal P}}'_2(\phi(0))\phi'(0)=\overline{{\cal P}}'_2(b)(a-bc)\]
(since $\overline{{\cal P}}_2$ is defined at $b$, then $\overline{{\cal P}}'_2$ is defined at $b$). We consider $\phi(t_2)\in {\Bbb K}(t_2)\setminus{\Bbb K}$ such that $$\overline{{\cal P}}'_2(b)\not=(0,0,0)\,\,\,\,\,\qquad \qquad \qquad \mbox{(III)}.$$

Observe that the rational function  $\phi$  satisfying (I), (II) and (III) exists since these three conditions determine a  non empty open subset $\Omega\subset \Bbb K$. Indeed: clearly, $\Omega$ is an open subset of $\Bbb K$. In addition,  if  $\Omega=\emptyset$, then $\overline{{\cal P}}'_2(t_2)=(0,0,0)$ which implies that  $\overline{{\cal P}}_2(t_2)\in{\Bbb K}^3$. This  is impossible since $\cal V$ is a surface parametrized by  $\overline{{\cal P}}(t_1, t_2)=\widehat{\cal P}_1(t_1)+\overline{{\cal P}}_2(t_2)$ and then $\overline{{\cal P}}_2(t_2)\not\in{\Bbb K}^3$. \para

Finally, we consider ${\cal P}(t_1,t_2)={\cal P}_1(t_1)+{\cal P}_2(t_2)$, where ${\cal P}_1(t_1)=\widehat{\cal P}_1(t_1)+\widehat{\cal P}_2(0)$, and ${\cal P}_2(t_2)=\widehat{\cal P}_2(t_2)-\widehat{\cal P}_2(0)$. Note that  ${\cal P}_2(t_2)\in {\Bbb K}(t_2)^3$ is such that  ${\cal P}_2(0)=(0,0,0)$,  ${\cal P}_2'(0)=\widehat{\cal P}'_2(0)=\overline{{\cal P}}'_2(b)(a-bc)\not=(0,0,0)$ (by (I), we have that $a-bc\not=0$). Furthermore $\cP_2$ is proper since $\widehat{\cal P}_2$ is proper. \qed

\para

\begin{remark}\label{t20} From the proof of Lemma \ref{L-prop1}, one deduces that we can have different  $t_2^0\in \Bbb{K}$ such that  ${\cal P}_2(t_2^0)=(0,0,0)$ and  ${\cal P}'_2(t_2^0)\not=(0,0,0)$ since the solution space $\Omega$ is a non empty open subset of  $\Bbb K$.  In addition, for every $t_2^0\in \Bbb{K}$ satisfying the above conditions, we obtain different proper parametrizations $\cP_1$ and $\cP_2$ solving the problem. In Section 4, we will take a deeper look at these $t_2^0$ to simplify our computation (see Theorem~\ref{T-MP1}).
\end{remark}


\section{Parameterizing the translational surface}

In this section, we provide  a necessary and sufficient condition for that an algebraic surface $\cal V$ is translational; i.e. $\cal V$ has a parametric representation of the form ${\cal P}(t_1,t_2)={\cal P}_1(t_1)+{\cal P}_2(t_2)$. The proof is constructive  and then, a method for computing $\cP$ is developed (see Theorem \ref{T-P1}, Corollary \ref{C-P1}, and statement 2 in Theorem \ref{T-NecSuf}).

\para

For this purpose, we assume that we are in the conditions stated in Section 1 and first, we prove Theorem \ref{T-P1} where we show a  necessary condition to that  $\cal V$  admits a parametrization of the form given in Eq. (\ref{eq-form}).  This condition provides a method for computing the parametrization  ${\cal P}_1$.

\para

\begin{theorem}\label{T-P1} Let  $\cal V$ be a translational surface. There exists  $(a_1, a_2, a_3)\in {\Bbb K}^3\setminus\{(0,0,0)\}$ such that
  ${\cal P}_1(t_1)\in {\Bbb K}(t_1)^3$ parametrizes properly a  space curve ${\cal C}_1\subset  {\Bbb K}^3$ defined by  the equations $f(\ox)=g(\ox)=0$, where $g(\ox):=a_1f_{x_1}(\ox)+a_2f_{x_2}(\ox)+a_3f_{x_3}(\ox)$.
\end{theorem}

\vspace*{2mm}

\noindent {\bf Proof.} Since $\cal V$  admits a   parametrization of the form given in Eq. (\ref{eq-form}), we have that $f({\cal P}_1(t_1)+{\cal P}_2(t_2))=0$, and in particular
\[0=f({\cal P}_1(t_1)+{\cal P}_2(0))=f({\cal P}_1(t_1))\]
(see Lemma \ref{L-prop1}). In addition, from $f({\cal P}_1(t_1)+{\cal P}_2(t_2))=0$, we also get that
$\nabla f({\cal P}_1(t_1)+{\cal P}_2(t_2))\,\cdot \,{\cal P}_2'(t_2)=0.$
Thus,   it holds that
\[0=\nabla f({\cal P}_1(t_1)+{\cal P}_2(0))\,\cdot \,{\cal P}_2'(0)=\nabla f({\cal P}_1(t_1))\,\cdot \,(a_1,a_2,a_3),\]
where ${\cal P}_2'(0)=(a_1,a_2,a_3)\in {\Bbb K}^3\setminus\{(0,0,0)\}$ (see Lemma \ref{L-prop1}).\para

Finally, we prove that the equations $f(\ox)=g(\ox)=0$ define a space curve. Indeed: since $\cal V$ is not a cylindrical surface, we get that $g(\ox):=a_1f_{x_1}(\ox)+a_2f_{x_2}(\ox)+a_3f_{x_3}(\ox)\not\in {\Bbb K}$ for every $(a_1, a_2, a_3)\in {\Bbb K}^3\setminus\{(0,0,0)\}$ (note that $g(\cP_1)=0$ and thus, if $g(\ox)=c\in {\Bbb K}$ then $c=0$ which would imply that $g=0$ and $\cal V$ is  a cylindrical surface). In addition, since $0<\deg(g)<\deg (f)$ and $f$ is irreducible, we get that $f, g$ are linearly independent. \qed

\para

\para

From Theorem \ref{T-P1}, one easily gets the following corollary where it is shown that in fact, there exists {\it infinitely many} vectors $(a_1,a_2,a_3)$  that can be used to compute a proper parametrization $\cP_1$  (as many as vectors given by ${\cal P}'_2(t_2^0)$).  Corollary \ref{C-P1} generalizes  Theorem \ref{T-P1} in the sense that   for each vector ${\cal P}'_2(t_2^0):=(a_1, a_2, a_3)$ we obtain a rational space curve ${\cal C}_1$ (this space curve will depend on the vector $(a_1,a_2,a_3)$ considered) such that a proper parametrization of ${\cal C}_1$ is given by ${\cal P}_1(t_1)+{\cal P}_2(t_2^0)$. Thus, $\overline{\cal P}_1+\overline{\cal P}_2$ is a parametrization of $\cal V$, where $\overline{\cal P}_1(t_1):={\cal P}_1(t_1)+{\cal P}_2(t_2^0)$ and $\overline{\cal P}_2(t_2):={\cal P}_2(t_2)-{\cal P}_2(t_2^0)$.

\para

\begin{corollary}\label{C-P1} Let  $\cal V$ be a translational surface.  For every   ${\cal P}'_2(t_2^0)=(a_1, a_2, a_3)\in {\Bbb K}^3\setminus\{(0,0,0)\}$ it holds that
  ${\cal P}_1(t_1)+{\cal P}_2(t_2^0)\in {\Bbb K}(t_1)^3$ parametrizes properly a  space curve  defined by  the equations $f(\ox)=g(\ox)=0$, where $g(\ox):=a_1f_{x_1}(\ox)+a_2f_{x_2}(\ox)+a_3f_{x_3}(\ox)$.  In addition,  $\overline{\cal P}_1+\overline{\cal P}_2$ is a parametrization of $\cal V$, where $\overline{\cal P}_1(t_1):={\cal P}_1(t_1)+{\cal P}_2(t_2^0)$ and $\overline{\cal P}_2(t_2):={\cal P}_2(t_2)-{\cal P}_2(t_2^0)$.
\end{corollary}

\para

\begin{remark}\label{pset-paramspace}
\begin{itemize}
\item[1.]  We prove Theorem~\ref{T-P1} and Corollary \ref{C-P1} in a constructive way, but we should remind that there may exists suitable $(a_1,a_2,a_3)$ not only from ${\cal P}_2'(t_2)$ (see Example \ref{E-5}). We prefer to give more discussions for $(a_1,a_2,a_3)$ from ${\cal P}_2'(t_2)$ since it is enough to lead a parametrization algorithm.
\item[2.] In order to compute a rational proper parametrization of a space curve ${\cal D}$, one may apply for instance the results in \cite{Chou1992} or in \cite{libro} (we remind that   any space curve  can be birationally projected onto a
plane curve).
\end{itemize}
\end{remark}

\para

In the following, we assume that the parametrization $\cP_1$ is computed (see Theorem \ref{T-P1} and Corollary \ref{C-P1}), and we show how to compute the parametrization $\cP_2$. For this purpose, we consider  $$h(\ox,t_1):=f({\cal P}_1(t_1)+(x_1,x_2,x_3))=\widetilde{h}(\ox)\Psi(\ox, t_1)p(t_1),$$
where $\Psi(\ox, t_1):=\widetilde{h}_0(\ox)+\widetilde{h}_1(\ox)t_1\cdots+\widetilde{h}_n(\ox)t_1^n\in {\Bbb K}[\ox, t_1]$,  $\gcd(\widetilde h_0,\ldots,\widetilde h_n)=1$, and $p(t_1)\in {\Bbb K}(t_1),\,\widetilde{h}(\ox)\in {\Bbb K}[\ox]$. We denote by   $V(\widetilde h_0,\ldots,\widetilde h_n)$    the variety generated by $\widetilde h_0,\ldots,\widetilde h_n$.

\para

In Theorem \ref{T-NecSuf}, we show that every rational space curve ${\cal D}\subset V(\widetilde h_0,\ldots,\widetilde h_n)$ provides a  parametrization ${\cal P}_2(t_2)\in {\Bbb K}(t_2)^3$ such that  $\cP(t_1, t_2)=\cP_1(t_1)+\cP_2(t_2)$ parametrizes $\cal V$. In fact, Theorem \ref{T-NecSuf} provides a necessary and sufficient condition for that algebraic surface has a parametrization of the form given in Eq. (\ref{eq-form}), and the proof is constructive (see statement $3$). More precisely, once a parametrization $\cP_1$ is computed (see  Corollary \ref{C-P1}), we apply Theorem \ref{T-NecSuf} to compute a parametrization $\cP_2$ (if it exists). If  $\cP_2$ does not exist, we can conclude that $\cal V$ is not a translational surface.

\para

In order to show Theorem \ref{T-NecSuf}, we first need to show the following   lemma where some properties of $\Psi(\ox, t_1), \widetilde{h}(\ox)$ and $p(t_1)$ are proved.

\para

\begin{lemma}\label{prop2} The following statements hold:
\begin{itemize}
\item[1.] $f({\cal P}_1(t_1)+\ox)\not=0$ (and then  $\widetilde{h}(\ox)\Psi(\ox, t_1)p(t_1)\not=0$). In addition, $p(t_1^0)\not=0$ for every $t_1^0\in {\Bbb K}$.
\item[2.] $\Psi(\ox, t_1)\in {\Bbb K}[\ox, t_1]\setminus {\Bbb K}[\ox]$.
\item[3.] $\Psi(\ox, t_1)\in {\Bbb K}[\ox, t_1]\setminus {\Bbb K}[t_1]$.
\end{itemize}
\end{lemma}

\noindent {\bf Proof.} \begin{itemize}
\item[1.] Let us assume that  $f({\cal P}_1(t_1)+\ox)=0$. Then,  $\nabla f({\cal P}_1(t_1^0)+\ox)\cdot(1,1,1)=0$ for every $t_1^0\in {\Bbb K}$ where $\cP_1$ is defined.  By applying the change of variable $\ox\rightarrow \ox-\cP_1(t_1^0)$, one   gets that  $\nabla f(\ox)\cdot(1,1,1)=0$ which contradicts  our assumption of $\cal V$ not being a cylindrical surface.\\
    In addition, we also have that $p(t_1^0)\not=0$ for every $t_1^0\in {\Bbb K}$. Otherwise, $f({\cal P}_1(t_1^0)+\ox)=0$ which implies that $\nabla f({\cal P}_1(t_1^0)+\ox)\cdot(1,1,1)=0$ which leads us to a contradiction reasoning as before.
\item[2.] Let us assume that  $\Psi(\ox, t_1) \in {\Bbb K}[\ox]$. That is, $f(\cP_1(t_1)+\ox)=\widetilde{h}(\ox)\Psi(\ox)p(t_1)$. Let $t_1^0\in {\Bbb K}$ such that  $\deg_{\ox}(f(\cP_1(t_1)+\ox))=\deg_{\ox}(f(\cP_1(t_1^0)+\ox))$. Then, we have that $$
p(t_1^0)f(\cP_1(t_1)+\ox)=p(t_1)f(\cP_1(t_1^0)+\ox)
$$
(note that from statement 1 above, we have that $p(t_1^0)\not=0$ for every $t_1^0\in {\Bbb K}$). Deriving w.r.t $t_1$, we get that
$$p(t_1^0)\nabla f({\cal P}_1(t_1)+\ox)\cdot\cP_1'(t_1)=f(\cP_1(t_1^0)+\ox) p'(t_1).$$
Since  $\deg_{\ox}(\nabla f({\cal P}_1(t_1)+\ox))<\deg_{\ox}(f(\cP_1(t_1)+\ox))=\deg_{\ox}(f(\cP_1(t_1^0)+\ox))$, we obtain that $p'(t_1)=0$ which implies that $p(t_1)=p(t_1^0)$. Thus,
$$f(\cP_1(t_1)+\ox)= f(\cP_1(t_1^0)+\ox)$$
and
$$ \nabla f({\cal P}_1(t_1^0)+\ox)\cdot\cP_1'(t_1)= \nabla f({\cal P}_1(t_1)+\ox)\cdot\cP_1'(t_1)=0.$$
By applying the change of variable $\ox\rightarrow \ox-\cP_1(t_1^0)$, one   gets that  $\nabla f(\ox)\cdot\cP_1'(t_1)=0$ and in particular,  $\nabla f(\ox)\cdot\cP_1'(t_1^0)=0$, for every $t_1^0\in {\Bbb K}$ where $\cP_1'(t_1^0)\in {\Bbb K}^3\setminus\{(0,0,0)\}$ exists. This contradicts  our assumption of $\cal V$ not being a cylindrical surface. Therefore, we conclude that $\Psi(\ox, t_1)\in {\Bbb K}[\ox, t_1]\setminus {\Bbb K}[\ox]$.
\item[3.]  In order to prove that $\Psi(\ox, t_1)\in {\Bbb K}[\ox, t_1]\setminus {\Bbb K}[t_1]$, one reasons similarly as in statement 2. \qed\end{itemize}




\para

\begin{theorem}\label{T-NecSuf}   The following statements are equivalent:
\begin{itemize}
\item[1.] $\cal V$ is a translational surface. 
\item[2.] There exists a proper parametrization ${\cal P}_2(t_2)\in {\Bbb K}(t_2)^3$ such that   $\cP(t_1, t_2)=\cP_1(t_1)+\cP_2(t_2)$ parametrizes  $\cal V$.
\item[3.] It holds that   ${\cal C}_2\subset V(\widetilde h_0,\ldots,\widetilde h_n)$. In addition, every rational space curve ${\cal D}\subset V(\widetilde h_0,\ldots,\widetilde h_n)$ provides a  parametrization ${\cal P}_2(t_2)\in {\Bbb K}(t_2)^3$ such that $\cP(t_1, t_2)=\cP_1(t_1)+\cP_2(t_2)$ parametrizes $\cal V$.
 \end{itemize}
\end{theorem}

\vspace*{2mm}

\noindent {\bf Proof.}\begin{itemize}
 \item[$1 \Leftrightarrow 2$]
 Clearly,   $\cal V$  admits  parametrization of the form given in Eq. (\ref{eq-form}) (i.e. $\cal V$ is a translational surface) if and only if there exists a parametrization ${\cal P}_2(t_2)\in {\Bbb K}(t_2)^3$ such that   $\cP(t_1, t_2)=\cP_1(t_1)+\cP_2(t_2)$ parametrizes  $\cal V$.
\item[$1 \Rightarrow 3$] If $\cal V$   admits a  parametrization of the form given in Eq. (\ref{eq-form}),  Theorem \ref{T-NecSuf} holds, and then $h({\cal P}(t_2), t_1)=0$. Thus,  $\Psi(\cP_2(t_2), t_1)=0$. Indeed: let us assume that $\Psi(\cP_2(t_2), t_1)\not=0$ and then  $\widetilde{h}(\cP_2(t_2))=0.$ Let us prove that this is impossible. For this purpose, we first consider   $t_1^0\in {\Bbb K}$ such that $\deg_{\ox}(f(\cP_1(t_1)+\ox))=\deg_{\ox}(f(\cP_1(t_1^0)+\ox))$. Now, taking into account that
\[
f(\cP_1(t_1)+\ox)=\widetilde{h}(\ox)\Psi(\ox, t_1)p(t_1),\qquad \mbox{(I)}
\]
we apply the change of variable $\ox\rightarrow \ox-\cP_1(t_1^0)$, and   we obtain that
$$f(\cP_1(t_1)+\ox-\cP_1(t_1^0))=\widetilde{h}(\ox-\cP_1(t_1^0))\Psi(\ox-\cP_1(t_1^0), t_1)p(t_1).$$
Then, for $t_1=t_1^0$ we get that $f(\ox)=\widetilde{h}(\ox-\cP_1(t_1^0))\Psi(\ox-\cP_1(t_1^0), t_1^0)p(t_1^0)$. Since $f(\ox)$ is irreducible and $\widetilde{h}(\ox)$ is not a constant (note that $\widetilde{h}(\cP_2(t_2))=0$ and thus, if $\widetilde{h}(\ox)=\beta \in {\Bbb K}$ then $\beta=0$; hence, $\widetilde{h}(\ox)=0$ which is impossible by statement 1 in Lemma \ref{prop2}), we get that   $\Psi(\ox-\cP_1(t_1^0), t_1^0)p(t_1^0)=\alpha\in {\Bbb K}\setminus\{0\}$ and then,
$f(\ox+\cP(t_1^0))=\alpha\widetilde{h}(\ox)$.  Substituting in (I), we obtain that
$$
f(\cP_1(t_1)+\ox)=f(\cP_1(t_1^0)+\ox)\Psi(\ox, t_1)p(t_1),\qquad \mbox{(II)}.
$$
Since $\deg_{\ox}(f(\cP_1(t_1)+\ox))=\deg_{\ox}(f(\cP_1(t_1^0)+\ox))$, one   gets that $\Psi(\ox, t_1)p(t_1)=c(t_1)$ with $c(t_1^0)=\alpha$. Deriving w.r.t $t_1$ in (II), we have that
$$\nabla f({\cal P}_1(t_1)+\ox)\cdot\cP_1'(t_1)=f(\cP_1(t_1^0)+\ox) c'(t_1).$$
But  $\deg_{\ox}(\nabla f({\cal P}_1(t_1)+\ox))<\deg_{\ox}(f(\cP_1(t_1^0)+\ox))$, hence $c'(t_1)=0$ which implies that $c(t_1)=c(t_1^0)=\alpha$ and
$$
f(\cP_1(t_1)+\ox)=f(\cP_1(t_1^0)+\ox)\Psi(\ox, t_1)p(t_1)=\alpha f(\cP_1(t_1^0)+\ox).
$$
Thus, from this equality we get that
$$
\alpha\nabla f({\cal P}_1(t_1^0)+\ox)\cdot\cP_1'(t_1)= \nabla f({\cal P}_1(t_1)+\ox)\cdot\cP_1'(t_1)=0.
$$
By applying the   change of variable $\ox\rightarrow \ox-\cP_1(t_1^0)$, we have that $\nabla f(\ox)\cdot\cP_1'(t_1)=0$ which can only happen for  $\cal V$ being a cylindrical surface, which contradicts  our assumption. Therefore, we conclude that  $\Psi(\cP_2(t_2), t_1)=0$.

\para

 Finally, since $\Psi(\cP_2(t_2), t_1)=0$ and  ${\cal P}_2(t_2)\in {\Bbb K}(t_2)^3\setminus {\Bbb K}^3$ does not depend on $t_1$, we get that $\widetilde h_i({\cal P}_2)=0,\,0\leq i\leq n$ which implies that ${\cal C}_2\subset V(\widetilde h_0,\ldots,\widetilde h_n)$.
\item[$3 \Leftarrow 1$] Let ${\cal M}(t_2)\in {\Bbb K}(t_2)^3$ be a parametrization of a rational space curve ${\cal D}\subset V(\widetilde h_0,\ldots,\widetilde h_n)$. Thus, $\widetilde h_i({\cal M})=0,\,i\in \{0,\ldots, n\}$. Then,  $h({\cal M}(t_2), t_1)=0$ and hence $f({\cal P})=0,$ where ${\cal P}(t_1,t_2)=\cP_1(t_1)+{\cal M}(t_2)$. Furthermore, according to Lemma~\ref{prop2} (see statements 2 and 3), $V(\widetilde h_0,\ldots,\widetilde h_n)$ defines an algebraic set independent to $t_1$, which means that ${\cal P}(t_1,t_2)$ defines a surface.
\qed
\end{itemize}

\noindent
In the following, we illustrate Theorems \ref{T-P1} and \ref{T-NecSuf} with an example.

\para

\begin{example}\label{E-1} Consider the surface $\cal V$ defined implicitly by the polynomial $f(\ox)=x_3+5x_1^2-6x_1x_2+2x_2^2\in {\Bbb C}[\ox].$
One may check that $\cal V$ satisfy the assumptions introduced in Section 2. Thus, we first determine a rational proper parametrization of the space curve ${\cal C}_1$ defined by the polynomials
\[f(\ox),\quad \mbox{and}\quad g(\ox):=a_1f_{x_1}(\ox)+a_2f_{x_2}(\ox)+a_3f_{x_3}(\ox) \]
(see Theorem \ref{T-P1}).  For this purpose, taking into account statement 2 in Remark \ref{pset-paramspace},  we compute a proper parametrization of the plane curve defined by
\[\res_{x_3}(f(\ox), g(\ox))=10a_1 x_1-6a_1 x_2-6 a_2 x_1+4 a_2 x_2+a_3.\]
We get that
${\cal Q}(t_1)=\left(t_1, {-10 a_1 t_1+6 a_2 t_1-a_3}/{2(-3 a_1+2 a_2)}\right)$
is a proper parametrization of the plane curve. Thus, a proper parametrization of ${\cal C}_1$ is given by
\[{\cal P}_1(t_1)=\left(t_1, \frac{-10 a_1 t_1+6 a_2 t_1-a_3}{2(-3 a_1+2 a_2)}, p(t_1)\right),\]
where $N(p(t_1))=0$, and
$N(x_3):=\gcd(f({\cal Q}(t_1),x_3), g({\cal Q}(t_1),x_3))$. We get
\[{\cal P}_1(t_1)=\left(t_1, \frac{-10 a_1 t_1+6 a_2 t_1-a_3}{2(-3 a_1+2 a_2)},  \frac{-12 t_1^2 a_1 a_2+4 t_1^2 a_2^2+2 t_1 a_3 a_1+10 t_1^2 a_1^2+a_3^2}{-2 (9 a_1^2-12 a_1 a_2+4 a_2^2)}\right).\]

\noindent
Now, we compute a rational proper parametrization of the space curve ${\cal C}_2$. For this purpose, we determine $\Psi(\ox, t_1)$, where
 $h(\ox,t_1):=f({\cal P}_1(t_1)+(x_1,x_2,x_3))=\widetilde{h}(\ox)\Psi(\ox, t_1)p(t_1).$  We get that\\

 \noindent
 $\Psi(\ox, t_1)=\widetilde{h}_0(\ox)+\widetilde{h}_1(\ox)t_1=3 a_1 x_3+6 a_1 x_2^2-18 a_1 x_1 x_2+15 a_1 x_1^2+12 x_1 x_2 a_2-3 x_1 a_3+2 x_2 a_3-2 x_3 a_2-4 x_2^2 a_2-10 x_1^2 a_2+(2 a_1 x_2-2 x_1 a_2) t_1.$\\

 \noindent
 Thus,  we compute a proper parametrization of the space curve defined by $V(\widetilde{h}_0, \widetilde{h}_1)$ (see Theorem \ref{T-NecSuf}).  Reasoning as above, we get that
\[{\cal P}_2(t_2)=\left(t_2, \frac{t_2 a_2}{a_1}, \frac{-t_2 (2 t_2 a_2^2+5 a_1^2 t_2-6 a_1 t_2 a_2-a_1 a_3)}{a_1^2}\right).\]
is a rational proper parametrization of  ${\cal C}_2$.\\

\noindent
Finally, we obtain that a parametrization of $\cal V$ is given by ${\cal P}(t_1, t_2)={\cal P}_1(t_1)+{\cal P}_2(t_2)=(p_1(t_1, t_2)/q_1(t_1, t_2), p_2(t_1, t_2)/q_2(t_1, t_2), p_3(t_1, t_2)/q_3(t_1, t_2)),$
\[p_1=t_1+t_2,\quad q_1=1,\quad p_2=6 a_1  t_2 a_2-4  t_2 a_2^2+10 a_1^2  t_1-6 a_1 a_2  t_1+a_1 a_3,\quad q_2=2 a_1 (3 a_1-2 a_2),\]

\noindent
$p_3=220  t_2^2 a_1^2 a_2^2-96  t_2^2 a_1 a_2^3+16  t_2^2 a_2^4+90  t_2^2 a_1^4-228  t_2^2 a_1^3 a_2-18  t_2 a_1^3 a_3+24  t_2 a_1^2 a_3 a_2-8  t_2 a_1 a_3 a_2^2-12  t_1^2 a_1^3 a_2+4  t_1^2 a_1^2 a_2^2+2  t_1 a_3 a_1^3+10  t_1^2 a_1^4+a_3^2 a_1^2,\quad q_3=-2a_1^2 (9 a_1^2-12 a_1 a_2+4 a_2^2).$

      \end{example}

\para

\para

In Example \ref{E-1}, we obtain a parametrization   $\cP(t_1, t_2)=\cP_1(t_1)+\cP_2(t_2)$ of the surface $\cal V$ that depends on the vector $(a_1, a_2, a_3)$ introduced in Theorem \ref{T-P1}. Observe that for a particular  value  of  this vector (such that the polynomials $f, g$ define a space curve), we obtain a particular parametrization of the form given in Eq. (\ref{eq-form}).  In Section 4, we deal with this question and   we show how the computation of $\cP$ can be (computationally speaking) improved by choosing any particular value for the vector $(a_1, a_2, a_3)$ (compare with Corollary \ref{C-P1}).  In particular, in Section 4,  two theorems are proved. The first one, Theorem \ref{T-P2}, shows how the computation of $\cP_2$ can be improved by considering a new (simpler) variety that generates the same space curves than  $V(\widetilde{h}_0,\ldots,\widetilde{h}_n)$.  The second theorem, Theorem \ref{T-MP1}, allow us to consider particular values for the parameters $a_1, a_2, a_3$ which makes easier the computation of the parametrizations $\cP_i,\,i=1,2$.


\section{Practical computation of the parametrization}

In Theorem \ref{T-NecSuf} (Section 3), we show how once the parametrization $\cP_1$ is computed (see Theorem \ref{T-P1} and Corollary \ref{C-P1}), one may determine  the parametrization $\cP_2$ (and then, a parametrization $\cP=\cP_1+\cP_2$ of the given surface $\cal V$) by computing a proper parametrization of a space curve ${\cal D} \subset V(\widetilde{h}_0,\ldots,\widetilde{h}_n)$, where    $$h(\ox,t_1):=f({\cal P}_1(t_1)+(x_1,x_2,x_3))=\widetilde{h}(\ox)\Psi(\ox, t_1)p(t_1),$$
 $\Psi(\ox, t_1):=\widetilde{h}_0(\ox)+\widetilde{h}_1(\ox)t_1\cdots+\widetilde{h}_n(\ox)t_1^n$,  and $\gcd(\widetilde h_0,\ldots,\widetilde h_n)=1$.  In fact, Theorem \ref{T-NecSuf} provides a necessary and sufficient condition for that algebraic surface having a parametrization of the form given in Eq. (\ref{eq-form}). More precisely, once a parametrization $\cP_1$ is computed (see  Corollary \ref{C-P1}), we apply Theorem \ref{T-NecSuf} to compute a parametrization $\cP_2$ (if it exists) by applying   statement $3$ in Theorem \ref{T-NecSuf}. If  $\cP_2$ does not exist, we can conclude that $\cal V$ is not a translational surface.

 \para

 In the following theorem, we show that the computation of $\cP_2$ can be improved in the sense that we do not need to compute explicitly $V(\widetilde{h}_0,\ldots,\widetilde{h}_n)$ but a simpler variety that generates the same space curves than $V(\widetilde{h}_0,\ldots,\widetilde{h}_n)$. More precisely, we prove that for {\it ``almost all'' }   pair of values    $s_1, s_2\in {\Bbb K}$, it holds that any rational space curve ${\cal D} \subset V(\widetilde{h}_0,\ldots,\widetilde{h}_n)$ can be defined by  the polynomials   $\Psi(\ox, s_i)\in {\Bbb K}[\ox],\,i=1,2.$
\para

\begin{theorem}\label{T-P2} Let ${\cal D}$ be a rational space curve such that ${\cal D} \subset V(\widetilde{h}_0,\ldots,\widetilde{h}_n)$.  There exists a non empty open subset $\Sigma\subset {\Bbb K}^2$ such that for every $(s_1, s_2)\in \Sigma$,  the space curve $\cal D$ is defined by  the polynomials  $g_i(\ox)\in {\Bbb K}[\ox]$, where $g_i(\ox):=\Psi(\ox, s_i)$ for $i=1,2.$

\end{theorem}

\vspace*{2mm}

\noindent {\bf Proof.} Let $\cP_2(t_2)=(p_{21}(t_2), p_{22}(t_2), p_{23}(t_2))$ be a proper parametrization of  ${\cal D}$, and let $R(x_1,x_2,y_1,y_2):=\res_{x_3}(\Psi(\ox, y_1), \Psi(\ox, y_2))$, where $y_1, y_2$ are new variables. Since $\Psi(\cP_2(t_2), t_1)=0$ for every $t_1$, we get that $R(\cQ(t_2),y_1,y_2)=0$ for every $y_1, y_2$, where $\cQ(t_2):=(p_{21}(t_2), p_{22}(t_2))$. Let
$$\Sigma:=\{(s_1, s_2)\in {\Bbb K}^2\,|\,R(x_1,x_2,s_1,s_2)\ell(x_1, x_2, s_1, s_2)\not=0\}\subset {\Bbb K}^2,$$  where $\ell\in {\Bbb K}[x_1, x_2, s_1, s_2]$ denotes the leader coefficient of $\Psi$ w.r.t $x_3$. Clearly, $\Sigma$ is a non empty open subset of ${\Bbb K}^2$. Now, we apply  the properties of specialization of resultants (see Lemma 4.3.1, p. 96 in \cite{win}), and we get that
\[R(x_1,x_2,s_1,s_2)=(\ell(x_1, x_2, s_1, s_2))^n\res_{x_3}(\Psi(\ox, s_1), \Psi(\ox, s_2)),\, n\in {\Bbb N}.\]
 Therefore, since $\Psi(\cP_2(t_2), s_i)=0,\,i=1,2$, and $\Psi(\ox, s_1), \Psi(\ox, s_2)$ are linearly independent (note that $R(x_1,x_2,s_1,s_2)\not=0$), we deduce that $R(\cQ(t_2),s_1,s_2)=0$ for every $(s_1, s_2)\in \Sigma$, and in particular $\cal D$ is defined by  the polynomials  $\Psi(\ox, s_i)\in {\Bbb K}[\ox],\,i=1,2$.
\qed

\para

\begin{remark}\label{R-P2}  Let $\cP_1(t_1)=(p_{11}(t_1), p_{12}(t_1), p_{13}(t_1))\in {\Bbb K}(t_1)^3,\,p_{1j}=\frac{p_{1j1}}{p_{1j2}},\,j=1,2,3$. Since $f(\cP_1(t_1)+\ox)=\widetilde{h}(\ox)\Psi(\ox, t_1)p(t_1)$,
it holds that, up to constants in ${\Bbb K}\setminus\{0\}$,  $$\Psi(\ox, s_i)=f(\cP_1(s_i)+\ox)/G(\ox),\,\,\,i=1,2,$$ where $G(\ox):=\gcd(f(\cP_1(s_1)+\ox), f(\cP_1(s_2)+\ox))$, and   $s_i\in {\Bbb K}$ is such that $\cP_1$ is defined (note that by statement 1 in by Lemma \ref{prop2}, we have that  $p(s_i)\not=0$). Therefore, for every $(s_1, s_2)\in \Omega$, where
$$\Omega:=\Sigma \cap \{(s_1, s_2)\in {\Bbb K}^2\,|\,p_1(s_1)p_1(s_2)\not=0\}\subset {\Bbb K}^2,\,\,\mbox{and $p_1:=\lcm(p_{112},p_{122},p_{132})$}$$
the space curve $\cal D$ is defined by  the polynomials   $g_i(\ox):=f(\cP_1(s_i)+\ox)/G(\ox),\,\,i=1,2.$ Note that $\Omega$ is again a non empty open subset of ${\Bbb K}^2$.
\end{remark}

\para

In the following example, we consider the surface $\cal V$ introduced in Example \ref{E-1}, and we show how the computation of $\cP_2$ can be improved by applying Theorem \ref{T-P2} or Remark \ref{R-P2}. In addition, we motivate the result that we will prove in Theorem \ref{T-MP1}, and we do not consider a generic vector $(a_1, a_2, a_3)$ (see Example \ref{E-1}). Instead, we take a particular value for this vector ($a_1=a_2=a_3=1$), and we show how a parametrization of the form ${\cal P}(t_1, t_2)={\cal P}_1(t_1)+{\cal P}_2(t_2)$ is provided for the input surface $\cal V$.\para

\begin{example}\label{E-2} Consider the surface $\cal V$ introduced in Example \ref{E-1}. $\cal V$ is implicitly defined by the polynomial $f(\ox)=x_3+5x_1^2-6x_1x_2+2x_2^2\in {\Bbb C}[\ox].$
\para

\noindent
First, we compute a parametrization of   the space curve ${\cal C}_1$ defined by the polynomials
\[f(\ox),\quad \mbox{and}\quad g(\ox):=a_1f_{x_1}(\ox)+a_2f_{x_2}(\ox)+a_3f_{x_3}(\ox) \]
(see Theorem \ref{T-P1}). We take the particular values  $a_1=a_2=a_3=1$, and we get that \[{\cal P}_1(t_1)=\left(t_1, \frac{1}{2}(4 t_1+1), -\frac{1}{2}(2 t_1^2 +2 t_1+1)\right) \] is a proper parametrization of ${\cal C}_1$.


\para

\noindent
Now, we apply  Remark \ref{R-P2}, and we compute a  parametrization of the space curve ${\cal C}_2$ defined by the polynomials $g_1(\ox)=f({\cal P}_1(1)+(x_1,x_2,x_3))$ and $g_2(\ox)=f({\cal P}_1(-3)+(x_1,x_2,x_3))$ (note that $G(\ox)=\gcd(f({\cal P}_1(1)+(x_1,x_2,x_3)), f({\cal P}_1(-3)+(x_1,x_2,x_3)))=1$).  We get that
\[{\cal P}_2(t_2)=\left(t_2, t_2, -t_2 ( t_2-1)\right)\]
is a rational proper parametrization of  ${\cal C}_2$.\para

\noindent
Finally, we obtain that a parametrization of $\cal V$ is given by
\[{\cal P}(t_1, t_2)={\cal P}_1(t_1)+{\cal P}_2(t_2)=\left(t_1+t_2, \frac{1}{2}(4 t_1+1)+t_2, -\frac{1}{2}(2 t_1^2 +2 t_1+1)-t_2 ( t_2-1)\right).\]
      \end{example}

      \para

In Example \ref{E-2}, we compute $\cP_1$ and $\cP_2$  for a particular value  $(a_1, a_2, a_3)\in {\Bbb K}^3$ (compare with Example \ref{E-1}). In fact,  we may take any value for this vector except  $a_1=a_2=0$. In this case,  $g(\ox):=a_1f_{x_1}(\ox)+a_2f_{x_2}(\ox)+a_3f_{x_3}(\ox)=a_3$ and thus, $f=g=0$ do not define a space curve. This remark gives us the idea that the computation of $\cP_1$ can be simplified since for particular values of the vector $(a_1, a_2, a_3)$, we can obtain a parametrization $\cP_1$. In order to show this property, we fist prove a technical lemma.

      \para

\begin{lemma}\label{L-MP1}  Let $\cal W$ be the surface defined by the parametrization
${\cal Q}(t_1, t_2)=(t_1 q_{1}(t_2), t_1 q_{2}(t_2), t_1 q_{3}(t_2))\in {\Bbb K}(t_1, t_2)^3.$ It holds that $(a_1, a_2, a_3)\in {\cal W}$, where $(a_1, a_2, a_3)\in {\Bbb K}^3\setminus\{(0,0,0)\}$,\, $a_1a_2a_3=0$ and $a_i=1$ for  $i\in \{1,2,3\}$.
\end{lemma}

\vspace*{2mm}

\noindent {\bf Proof.}
 According to the construction, ${\cal Q}(t_1,t_2)$ is a rational conical surface, and then the implicit polynomial defining $\cal W$, $f_{\cal W}(x_1,x_2,x_3)$, is homogenous w.r.t $x_1,x_2,x_3$.
Then without loss of generality, setting $x_1=0$, we get  $\overline{f}(x_2,x_3):=f_{\cal W}(0,x_2,x_3)=0$ is homogenous in $x_2$ and $x_3$. Then, if $(a,b)\in {\Bbb K}^2\setminus\{(0,0)\}$  is such that $\overline{f}(a,b)=0$, we get that $f_{\cal W}(0,1,b/a)=\overline{f}(1,b/a)=0$ (if $a\not=0$) or $f_{\cal W}(0,a/b,1)=\overline{f}(a/b,1)=0$ (if $b\not=0$). Thus, $(0,1,b/a)\in {\cal W}$ or $(0,a/b, 1)\in {\cal W}$. Similar discussions in the cases of setting $x_2=0$ or $x_3=0$.
  \qed

\para

In the following, we consider the surface $\cal W$ defined by the parametrization ${\cal Q}(t_1, t_2)=(t_1 q_{1}(t_2), t_1 q_{2}(t_2), t_1 q_{3}(t_2))$, where  $\cP_2'(t_2):=(q_1(t_2), q_2(t_2), q_3(t_2))$. Observe that $\cal Q$ parametrizes a surface since  its jacobian has rank 2; otherwise, $q_iq'_j=q'_iq_i$ which implies that $(q_i/q_j)'=0$ and then $q_i=c_j q_j$ for  $i,j\in \{1,2,3\}$, $i\not=j$ and $c_j\in {\Bbb K}$. Thus,
${\cal P}_2(t_2)=(p(t_2), \alpha_1 p(t_2)+\beta_1, \alpha_2 p(t_2)+\beta_2),\,\, \alpha_i, \beta_i\in{\Bbb K},$ which is impossible since $\cP_2$ is not a line (see statement 4 in Proposition \ref{P-properties}).

\para

Under these conditions, Theorem \ref{T-MP1} shows that we can improve the computation of $\cP_1$ by taking a particular value for the vector $(a_1, a_2, a_3)$ introduced in Theorem \ref{T-P1} (compare Examples \ref{E-1} and \ref{E-2}). We illustrate this result with  Examples \ref{E-3} and \ref{E-4}.

  \para

\begin{theorem}\label{T-MP1} Let  $\cal V$ be a translational surface.  For every   $(a_1, a_2, a_3)\in {\cal W}$ it holds that there exists $t_2^0\in {\Bbb K}$ such that
  ${\cal P}_1(t_1)+{\cal P}_2(t_2^0)\in {\Bbb K}(t_1)^3$ parametrizes properly a  space curve  defined by  the equations $f(\ox)=g(\ox)=0$, where $g(\ox):=a_1f_{x_1}(\ox)+a_2f_{x_2}(\ox)+a_3f_{x_3}(\ox)$.  In addition,  $\overline{\cal P}_1+\overline{\cal P}_2$ is a parametrization of $\cal V$, where $\overline{\cal P}_1(t_1):={\cal P}_1(t_1)+{\cal P}_2(t_2^0)$ and $\overline{\cal P}_2(t_2):={\cal P}_2(t_2)-{\cal P}_2(t_2^0)$.
\end{theorem}

\vspace*{2mm}

\noindent {\bf Proof.} Let $p:=(a_1, a_2, a_3)\in {\cal W}$. In order to prove this theorem, we distinguish two different cases:
\begin{itemize}
\item[1.] There exists $(t_1^0, t_2^0)\in {\Bbb K}^2$ such that ${\cal Q}(t_1^0, t_2^0)=p$. We assume w.l.o.g. that $t_1^0=1$ (otherwise, we consider the reparametrization ${\cal Q}(t_1t_1^0, t_2)$). Thus, we have that $\cP_2'(t_2^0)=p$. Now, we apply Corollary \ref{C-P1} and we conclude that
  ${\cal P}_1(t_1)+{\cal P}_2(t_2^0)\in {\Bbb K}(t_1)^3$ parametrizes properly a  space curve  defined by  the equations $f(\ox)=g(\ox)=0$, where $g(\ox):=\nabla f(\ox)\cdot p$.  In addition,  $\overline{\cal P}_1+\overline{\cal P}_2$ is a parametrization of $\cal V$, where $\overline{\cal P}_1(t_1):={\cal P}_1(t_1)+{\cal P}_2(t_2^0)$ and $\overline{\cal P}_2(t_2):={\cal P}_2(t_2)-{\cal P}_2(t_2^0)$.
\item[2.] There does not exist $(t_1^0, t_2^0)\in {\Bbb K}^2$ such that ${\cal Q}(t_1^0, t_2^0)=p$.   Then, we  consider a reparametrization ${\cal Q}^*(t,s)$ such that ${\cal Q}^*(t^0, s^0)=p$, where ${\cal Q}^*(t,s)={\cal Q}(R(t,s))=R_1(t,s)\cP_2'(R_2(t,s))$, and $R(t,s):=(R_1(t,s),R_2(t,s))\in ({\Bbb K}(t,s) \setminus{\Bbb K})^2$. Note that ${\cal Q}^*$ is again a parametrization of the surface $\cal W$.   Under these conditions, since $f(\cP_1(t_1)+\cP_2(R_2(t,s)))=0$, we get that $$\nabla f(\cP_1(t_1)+\cP_2(R_2(t,s))) \, \cdot \, \cP_2'(R_2(t,s)) R_2'(t,s)=0.$$ Then, $\nabla f(\cP_1(t_1)+\cP_2(R_2(t,s)))\, \cdot \,\cP_2'(R_2(t,s))=0$, and thus $$\nabla f(\cP_1(t_1)+\cP_2(R_2(t,s))) \, \cdot \,   \cP_2'(R_2(t,s)) R_1(t,s)=0.$$ Hence, $\nabla f(\cP_1(t_1)+\cP_2(R_2(t,s))) \, \cdot \, {\cal Q}^*(t,s)=0$, and in particular $$\nabla f(\cP_1(t_1)+\cP_2(R_2(t^0,s^0))) \, \cdot \, {\cal Q}^*(t^0,s^0)=\nabla f(\cP_1(t_1)+\cP_2(R_2(t^0,s^0)))\, \cdot \,  p=0.$$ Therefore, $\cP_1(t_1)+\cP_2(R_2(t^0,s^0)) \in {\Bbb K}(t_1)^3$ parametrizes properly a  space curve  defined by  the equations $f(\ox)=g(\ox)=0$, where $g(\ox):=\nabla f(\ox)\cdot p$.  In addition,  $\overline{\cal P}_1+\overline{\cal P}_2$ is a parametrization of $\cal V$, where $\overline{\cal P}_1(t_1):={\cal P}_1(t_1)+\cP_2(R_2(t^0,s^0))$ and $\overline{\cal P}_2(t_2):={\cal P}_2(t_2)-\cP_2(R_2(t^0,s^0))$. \qed
\end{itemize}

\para

      \begin{remark} \label{R-MP1} \begin{itemize} \item[1.] Theorem \ref{T-MP1} provides a necessary (but not sufficient) condition for the computation of $\cP_1$.
      \item[2.] Taking into account Lemma \ref{L-MP1}, one may apply Theorem \ref{T-MP1} for some vector of the form  $(a_1, a_2, a_3)\in {\Bbb K}^3\setminus\{(0,0,0)\}$ with $a_1a_2a_3=0$ and $a_i=1$ for some $i\in \{1,2,3\}$.
      \end{itemize}
      \end{remark}

\para

\begin{example}\label{E-3} We consider the surface $\cal V$ over $\Bbb C$  implicitly defined by the polynomial\\

\noindent
$f(\ox)=2 x_1 x_3^3 x_2-2 x_1 x_3^2 x_2-3 x_1 x_3 x_2+4 x_1^2 x_3 x_2+10 x_3^2 x_2+5 x_1 x_2+2 x_1 x_3^5-x_3^3 x_2^2+x_3^3 x_1^3+2 x_3^5 x_2-x_1^3 x_2-4 x_1^2 x_3^4-x_3^7-6 x_2^2+x_1^3-x_1^2-8 x_3^4-2 x_3^5-15 x_3^2-x_3^3+x_2^3+12 x_2+5 x_1 x_3^2+9 x_3^3 x_1-4 x_1^2 x_3+2 x_3^3 x_2-2 x_1 x_2^2+x_1^2 x_2+6 x_3 x_1+4 x_3 x_2-x_3^3 x_1^2-11 x_3-2 x_1-3 x_1^2 x_3^2+x_3^4 x_2-2 x_3^2 x_2^2-9.$\\

\noindent
Now, let ${\cal C}_1$ be the space defined by the polynomials $f(\ox)$ and $g(\ox):=f_{x_3}(\ox)=$\\

\noindent
$-11+6 x_1 x_3^2 x_2+6 x_1-30 x_3+4 x_2+4 x_1^2 x_2-16 x_3^3 x_1^2+10 x_3^4 x_2+6 x_3^2 x_2-3 x_1 x_2+27 x_1 x_3^2-6 x_1^2 x_3+4 x_3^3 x_2+10 x_3 x_1+20 x_3 x_2-3 x_1^2 x_3^2-3 x_3^2 x_2^2+3 x_3^2 x_1^3+10 x_1 x_3^4-4 x_3 x_2^2-4 x_1^2-3 x_3^2-10 x_3^4-32 x_3^3-7 x_3^6-4 x_1 x_3 x_2.$\\

\noindent
(see Theorem  \ref{T-MP1} and statement $2$ in Remark \ref{R-MP1}). By applying statement 2 in Remark \ref{pset-paramspace} (see also Example \ref{E-1}), we compute a proper rational parametrization of ${\cal C}_1$. We get
\[{\cal P}_1(t_1)=\left(t_1, \frac{1+t_1^2}{t_1^2}, \frac{1}{t_1}\right).\]

\noindent
Now, let ${\cal C}_2$ be the space curve defined by the polynomials $g_1(\ox)=f({\cal P}_1(1)+(x_1,x_2,x_3))$ and $g_2(\ox)=f({\cal P}_1(-3)+(x_1,x_2,x_3))$ (note that $G(\ox)=\gcd(f({\cal P}_1(1)+(x_1,x_2,x_3)), f({\cal P}_1(-3)+(x_1,x_2,x_3)))=1$; see Theorem \ref{T-P2}, and Remark \ref{R-P2}).  We get that a rational proper parametrization of  ${\cal C}_2$ is given by
\[{\cal P}_2(t_2)=\left(t_2^2, t_2^3, t_2\right).\]

\noindent
Finally, we obtain that a parametrization of $\cal V$ is given by
\[{\cal P}(t_1, t_2)={\cal P}_1(t_1)+{\cal P}_2(t_2)=\left(t_1+t_2^2, \frac{1+t_1^2}{t_1^2}+t_2^3, \frac{1}{t_1}+t_2\right)\in {\Bbb C}(t_1,t_2)^3.\]
      \end{example}

\para

In the following, we provide an example that clarifies statement $2$ in the proof of Theorem \ref{T-MP1}.

\para

\begin{example}\label{E-4} Let $\cal V$ be the surface over $\Bbb C$ implicitly defined by the polynomial
\[f(\ox)=
  x_1^4-2\,x_3+7\,x_3x_1+2\,x_2^2-5\,x_2x_3+x_
3^2+2\,x_1^3-10\,x_1^2x_2-2\,x_3x_1
^2+7\,x_1x_2^2-x_2^3.\]
According to statement $2$ in Remark \ref{R-MP1}, one can select  $(a_1,a_2,a_3)=(1,0,0)$. Reasoning similarly as in Example~\ref{E-3}, we   get a proper parametrization of ${\cal C}_1$ given by
\begin{equation}\label{E-4eq1}
  {\cal P}_1(t_1)=(t_1,t_1,t_1^2),
\end{equation}
and a proper parametrization of ${\cal C}_2$  given by
\begin{equation}\label{E-4eq2}{\cal P}_2(t_2)=(t_2,t_2^2,t_2^3).\end{equation}
Then, ${\cal P}(t_1,t_2)={\cal P}_1(t_1)+{\cal P}_2(t_2)$ is a proper parametrization of the surface $\cal V$.  Observe that ${\cal Q}(t_1,t_2)=(t_1, 2t_1t_2,3t_1t_2^2)$ and ${\cal Q}(1,0)=(1,0,0).$ Furthermore,  ${\cal P}_2(0)=(0,0,0)$ and ${\cal P}'_2(0)=(1,0,0)$ (see Lemma \ref{L-prop1}).
 \\

 \para

Now, we consider a different vector $(a_1,a_2,a_3)=(0,0,1)$, and we get  the space curve defined by the polynomials $f(\ox)$ and $g(\ox):=f_{x_3}(\ox)=
-2+7\,x_1-5\,x_2+2\,x_3-2\,x_1^2$. Here, reasoning as above, we  get the proper parametrizations
$${\cal P}_1(t_1)=\left (t_1,t_1-\frac{1}{4},\frac{3}{8}-t_1+t_1^2\right ),\quad \mbox{and}\quad {\cal P}_2(t_2)=\left(t_2,t_2+t_2^2,\frac{3}{4}\,t_2+\frac{3}{2}\,t_2^2+t_2
^3\right).$$ Thus, ${\cal P}(t_1,t_2)={\cal P}_1(t_1)+{\cal P}_2(t_2)$ is a proper parametrization of  $\cal V$.  However, ${\cal Q}(t_1,t_2)=(t_1,t_1(1+2\,t_2), t_1(\frac{3}{4}+3\,t_2+3\,t_2^2))$ and  obviously, there does not exist $(t_1^0,t_2^0)$ such that ${\cal Q}(t_1^0,t_2^0)=(0,0,1).$ Note that ${\cal Q}$ parametrizes the surface $\cal W$ defined by the polynomial $4x_1x_3-3x_2$, and $(0,0,1)\in {\cal W}.$
 \para

Observe that if we set ${\cal Q}(1,-\frac{1}{2})=(1,0,0)$, we get that $(a_1,a_2,a_3)=(1,0,0)$, and in this case, we get the parametrizations given in  Eq.~(\ref{E-4eq1}) and Eq.~(\ref{E-4eq2}). In addition, taking into account Theorem \ref{T-MP1}, we   also get parametrizations by computing
\begin{equation}\label{E-4eq3}\overline{\cal P}_1(t_1)={\cal P}_1(t_1)+{\cal P}_2(t_2^0)=(t_1-\frac{1}{2}, (t_1-\frac{1}{2}), (t_1-\frac{1}{2})^2),\end{equation}
and
\begin{equation}\label{E-4eq4}\overline{\cal P}_2(t_2)={\cal P}_2(t_2)-{\cal P}_2(t_2^0)=\left(t_2+\frac{1}{2},(t_2+\frac{1}{2})^2, (t_2+\frac{1}{2})^3\right).\end{equation}
Now $\overline{\cal P}(t_1,t_2)=\overline{\cal P}_1(t_1)+\overline{\cal P}_2(t_2)$ is a parametrization of $\cal V$, and  $\overline{\cal P}_2(-\frac{1}{2})=(0,0,0)$ and $\overline{\cal P}'_2(-\frac{1}{2})=(1,0,0)$, which support Lemma~\ref{L-prop1} and Theorem~\ref{T-P1}.\para

Finally, one should note that Eq.~(\ref{E-4eq1}) and Eq.~(\ref{E-4eq3}) define the same space curve ${\cal C}_1$. Similarly,  Eq.~(\ref{E-4eq2}) and Eq.~(\ref{E-4eq4}) are proper parametrizations of the same space curve ${\cal C}_2$.
\end{example}

\para

Finally, we give an example to illustrate  statement 1 in Remark~\ref{pset-paramspace}. More precisely,  Theorem~\ref{T-P1} and Corollary \ref{C-P1} are proved in a constructive way, taking $(a_1,a_2,a_3)={\cal P}_2'(t_2^0)$, for any $t_2^0\in {\Bbb K}$. However,   there may exist suitable $(a_1,a_2,a_3)$ not only from ${\cal P}_2'(t_2)$.

\para

\begin{example}\label{E-5} Let $\cal V$ be the surface  of Example~\ref{E-4} implictly defined by the polynomial\\
\noindent $f(\ox)=
  x_1^4-2\,x_3+7\,x_3x_1+2\,x_2^2-5\,x_2x_3+x_
3^2+2\,x_1^3-10\,x_1^2x_2-2\,x_3x_1
^2+7\,x_1x_2^2-x_2^3.$
If we set $(a_1,a_2,a_3)=(1,1,1)$, we get the proper parametrizations of ${\cal C}_1$ and  ${\cal C}_2$ given by
$${\cal P}_1(t_1)=\left (t_1,t_1-\frac{1}{4},\frac{3}{8}-t_1+t_1^2\right ),\,\,\,\mbox{and}\,\,\,\,{\cal P}_2(t_2)=\left(t_2,t_2+t_2^2,\frac{3}{4}\,t_2+\frac{3}{2}\,t_2^2+t_2
^3\right).$$
 One can check that ${\cal P}(t_1,t_2)={\cal P}_1(t_1)+{\cal P}_2(t_2)$ is a parametrization of the surface $\cal V$.  However, ${\cal Q}(t_1,t_2)=(t_1,t_1(1+2\,t_2), t_1(\frac{3}{4}+3\,t_2+3\,t_2^2))$ parametrizes the surface $\cal W$ defined by the polynomial $4x_1x_3-3x_2$, and  $(1,1,1)\not\in {\cal W}$.

 \para

Nevertheless, one can find an alternative $(a_1,a_2,a_3)$ and then compute ${\cal P}_1(t_1)$ and ${\cal P}_2(t_2)$  of the forms given in Lemma~\ref{L-prop1} and Theorem~\ref{T-P1} (compare with Example~\ref{E-4}).
\end{example}

\end{document}